\documentclass[12pt,reqno]{amsart}

\newcommand{\R}{{\mathbb R}}

\newcommand{\e}{\epsilon}

\newcommand{\Om}{\Omega}

\newcommand{\p}{\partial}
\renewcommand{\phi}{\varphi}

\numberwithin{equation}{section}
\newtheorem{theorem}{Theorem}[section]
\newtheorem{defn}[theorem]{Definition}
\newtheorem{lemma}[theorem]{Lemma}
\newtheorem{remark}[theorem]{Remark}

\newtheorem{corollary}[theorem]{Corollary}

\begin{document}

\title[Effective macroscopic dynamics of SPDEs ]
{Effective macroscopic dynamics of stochastic partial differential
equations in perforated domains }

\author[W. Wang, D. Cao \& J. Duan   ]
{Wei  Wang,  Daomin Cao  \&  Jinqiao Duan }

\address[W.~Wang ]
{Institute of Applied Mathematics\\ Chinese Academy of Sciences\\
Beijing, 100080, China } \email[W.~Wang]{wangwei@amss.ac.cn}

\address[D.~Cao ]
{Institute of Applied Mathematics\\ Chinese Academy of Sciences\\
Beijing, 100080, China } \email[D.~Cao]{dmcao@amt.ac.cn}

\address[J.~Duan]
{Department of Applied Mathematics\\
Illinois Institute of Technology\\
Chicago, IL 60616, USA} \email[J.~Duan]{duan@iit.edu}

\date{December 31, 2005 submitted; June 20, 2006  Revised; and August 30, 2006 accepted.}

\thanks{This work was partly supported by the NSF Grants
DMS-0209326 \& DMS-0542450 and the Outstanding Overseas Chinese
Scholars Fund of the Chinese Academy of Sciences.}


\begin{abstract}
An effective macroscopic model for a   stochastic microscopic
system is derived. The original microscopic system is modeled by a
stochastic partial differential equation   defined on a domain
perforated with small holes or heterogeneities. The homogenized
effective model is still a stochastic partial differential
equation but  defined on a unified   domain  without holes.  The
solutions of the microscopic model is shown to converge to those
of the effective macroscopic model in probability distribution, as
the size of holes diminishes to zero. Moreover, the long time
effectivity of the    macroscopic system in the sense of
\emph{convergence in  probability distribution}, and the
effectivity of the macroscopic system in the sense of
\emph{convergence in energy} are also proved.

\end{abstract}

\maketitle

\section{Introduction}\label{s1}

 In recent years there has been   explosive growth of activities in
   multiscale modeling of complex phenomena  in many areas
 including material science, climate dynamics, chemistry and biology \cite{E00, TM05}.
 Stochastic partial differential equations (SPDEs or stochastic PDEs) --- evolutionary
 equations containing   noises ---   arise naturally as
 mathematical models of    multiscale systems under random
 influences.  In fact the need to include stochastic
effects in mathematical modelling of realistic  physical behavior
has become widely recognized in, for example, condensed matter
physics, climate and geophysical sciences, and materials sciences.
But implementing this idea poses some challenges both in theory and
for computation \cite{Imkeller, WaymireDuan}.

 This paper is devoted to the effective macroscopic dynamics of
 microscopic systems modeled by
 parabolic SPDEs in perforated media which exhibit  small-scale
heterogeneities. One example of such microscopic systems of
interest is   composite materials with microscopic heterogeneities
under the impact of external random fluctuations.

\bigskip

The heterogeneity scale is taken to be much smaller than the
macroscopic scale, which is equivalent, here, to assuming that the
heterogeneities are evenly distributed. From a mathematical point
of view, one can  assume that microscopic heterogeneities (holes)
are periodically placed in the media. This periodicity can be
represented by a small positive parameter $\e$ (i.e., the period).
In fact we work on the space-time cylinder  $D_\e\times(0,T)$,
with $T>0$, and $D_\e$ being the spatial domain obtained by
removing a number $N_\e$ of holes, of size $\e$, periodically
distributed, from a fixed domain $D$. When taking $\e\rightarrow 0
$, the holes covering $D$ are smaller and smaller and their
numbers goes to $\infty$. This signifies that the heterogeneities
are finer and finer.

There are lots of work on the homogenization problem for the
deterministic systems defined in such perforated domain or other
heterogeneous media, see for example \cite{BOFM92, NR1, NR2,
SK02,TTM02} for heat transfer in a composite material,
\cite{BOFM92,CD89,CDMZ91} for the wave propagation in a composite
material and \cite{LiM05, MP99} for the fluid flow in a porous
media. For an  introduction see \cite{CD99, JKO94, SP80}.

Recently there are also   works on   homogenization in the random
context; see \cite{KP02, MM86, PP03} for general random
coefficients, and \cite{BM98, ZV93, ZV94} for randomly perforated
domains. And also see a survey book about the homogenization
results in a random context \cite{JKO94}. A basic assumption in
these texts is the ergodic hypotheses on the coefficients for the
passing of the limit of $\e \to 0$. Note that the microscopic
models in these works are partial differential equations with
random coefficients, so-called random partial differential
equations (random PDEs), instead of stochastic PDEs.

\bigskip

In the present paper, the microscopic model is a SPDE defined in a
perforated domain.  Homogenization techniques are   employed to
derive an effective, simplified,  macroscopic model.
 Homogenization is a formal mathematic procedure  for deriving
 macroscopic models from microscopic systems.
 It has been applied to a   variety of problems including composite
materials modelling, porous media and climate modelling; see
\cite{CD99,CK97,JKO94,SP80}. Homogenization provides effective
macroscopic behavior of the system with microscopical
heterogeneities for which  direct numerical simulations are
usually too expensive.

We consider a spatially extended system where stochastic effects are
taken into account in the model equation, defined on a deterministic
domain but perforated with small scale holes. Specifically, we study
a class of stochastic partial differential equations driven by white
noise on a perforated domain   in the following  form
\begin{equation*} \label{a}
du_\e(t)=(\mathcal{A}_\e u_\e+F_\e(x, t))dt+G_\e(x, t) dW(t), \;\;
0< t < T , \e >0,
\end{equation*}
which will be described  in detail  in the next section.
 For the general theory of SPDEs we refer to \cite{PZ92}.
 The goal here is to derive the homogenized equation (effective
 equation), which is still a stochastic partial differential equation,  for
(\ref{a}) by the homogenization techniques in the sense of
\emph{probability}.

Homogenization theory has been   developed for deterministic
systems, and compactness discussion for the solutions
$\{u_\e\}_\e$ in some function space is a key step in the
homogenized approach \cite{CD99}. However, due to the appearance
of the stochastic term in the above microscopic system considered
in this paper, such compact result does not hold  for this
stochastic system. Fortunately the compactness  in the sense of
probability, that is, the tightness of the distributions for
$\{u_\e\}$, still holds. So one appropriate way is to homogenize
the stochastic system in the sense of probability. The goal in
this paper is to derive an effective macroscopic equation for the
above microscopic system, by homogenization in the sense of
probability.  It is shown that the solution $u_\e$ of the
microscopic or heterogeneous system converges to that of the
macroscopic or homogenized system as $\e\downarrow 0$ in
probability. That also implies the distributions of $\{u_\e \}_\e$
weakly converge, in some appropriate space, to the distribution of
a stochastic process which solves the macroscopic effective
equation. Moreover, the long time effectivity  of the homogenized
macroscopic system is demonstrated, that is, the solution
$u_\e(t)$ is shown to converge to the stationary solution of the
homogenized equation as $t\rightarrow \infty$ and $\e\downarrow 0$
in the sense of probability distribution. Furthermore, the
effectivity  of the
  macroscopic system in the sense of convergence in energy is also
shown.

In our approach, one difficulty is that the spatial domain is
changing as $\e\rightarrow 0$. To overcome this we use the
extension operator introduced in \cite{CD89} and introduce a new
probability space depending on a parameter in which the solution
is uniformly bounded. One novelty here is that the original
microscopic model is a stochastic PDE, instead of a random PDE as
studied by others, e.g., \cite{KP02, PP03,Souganidis}.


\medskip

This paper is organized as follows. The   problem formulation is
stated in  \S \ref{s2}. Section \ref{s3} is devoted to basic
properties of the microscopic system. The   effective macroscopic
equation is derived in \S \ref{s4}.  The long time effectivity of
the homogenized macroscopic system is considered in \S \ref{s6}.
Finally, the effectivity of the macroscopic system in the sense of
convergence in energy is shown in \S \ref{s5}. Moreover, in the
Appendix we present the explicit expression of
the  homogenization matrix.\\


\section{Problem formulation }\label{s2}

Let $D$ be an open bounded set in $\R^n$, $n\geq 2$, with smooth
boundary $\p D$ and $\e>0$ is a small parameter. Let
$Y=[0,l_1)\times [0,l_2)\times\cdots\times [0,l_n)$ be a
representative (cubic) cell in $\R^n$ and $S$ an open subset of
$Y$ with smooth boundary $\p S$, such that $\overline{S}\subset
Y$. Write $l=(l_1, l_2, \cdots, l_n)$. Define $\e S=\{\e y: y\in
S\}$. Denote by $S_{\e,k}$ the translated image of $\e S$ by $kl$,
$k\in Z^n$, $kl=(k_1l_1, k_2l_2,\cdots, k_nl_n)$.  And let $S_\e$
be the set of all the holes contained in $D$ and $D_\e=D\backslash
S_\e$. Then $D_\e$ is a periodically perforated domain with holes
of the same size as period $\e$. We   assume that the holes do not
  intersect with the boundary $\p D$, which implies that $\p
D_\e=\p D \cup \p S_\e$. See Fig. 1 for the case $n=2$. This
assumption is for avoiding technicalities and the results of our
paper will remain valid without this assumption; see \cite{AMN93}.

In the sequel we use the notations
\begin{equation*}
Y^*=Y\backslash\overline{S}, \;\;\vartheta=\frac{|Y^*|}{|Y|}
\end{equation*}
with $|Y|$ and $|Y^*|$   the  Lebesgue measure of $Y$ and $Y^*$
respectively. And denote by $\tilde{v}$ the zero extension to the
whole $D$ for any function defined on $D_\e$:

$$ \tilde{v}= \left\{
  \begin{array}{c l}
     v & \mbox{on  $D_\e$},\\
     0 & \mbox{on $S_\e$}.
  \end{array}
\right.  $$

\vskip 1cm
\begin{center}
 \setlength{\unitlength}{.35cm}
\begin{picture}(35,13)


  \qbezier(0,2)(0, 1)(1,1)
  \qbezier(0,10)(0,1)(0,2)
  \qbezier(0,10)(0,11)(1,11)
 \qbezier(1,11)(15,11)(15,11)
 \qbezier(15,11)(16,11)(16,10)
 \qbezier(16,10)(16,2)(16,2)
 \qbezier(16,2)(16,1)(15,1)
 \qbezier(15,1)(1,1)(1,1)


\qbezier(.5,1.8)(.5,0.8)(2.5,2.3)
 \qbezier(.5,1.8)(.5,2.6)(.5,2.6)
 \qbezier(.8,2.9)(.5,2.9)(2.5,2.9)
 \qbezier(.5,2.6)(.5,2.9)(.8,2.9)
 \qbezier(2.5,2.9)(3.5,2.9)(2.5,2.3)

 \qbezier(3.5,1.8)(3.5,0.8)(5.5,2.3)
 \qbezier(3.5,1.8)(3.5,2.6)(3.5,2.6)
 \qbezier(3.8,2.9)(3.5,2.9)(5.5,2.9)
 \qbezier(3.5,2.6)(3.5,2.9)(3.8,2.9)
 \qbezier(5.5,2.9)(6.5,2.9)(5.5,2.3)

 \qbezier(6.5,1.8)(6.5,0.8)(8.5,2.3)
 \qbezier(6.5,1.8)(6.5,2.6)(6.5,2.6)
 \qbezier(6.8,2.9)(6.5,2.9)(8.5,2.9)
 \qbezier(6.5,2.6)(6.5,2.9)(6.8,2.9)
 \qbezier(8.5,2.9)(9.5,2.9)(8.5,2.3)

 \qbezier(9.5,1.8)(9.5,0.8)(11.5,2.3)
 \qbezier(9.5,1.8)(9.5,2.6)(9.5,2.6)
 \qbezier(9.8,2.9)(9.5,2.9)(11.5,2.9)
 \qbezier(9.5,2.6)(9.5,2.9)(9.8,2.9)
 \qbezier(11.5,2.9)(12.5,2.9)(11.5,2.3)

 \qbezier(12.5,1.8)(12.5,0.8)(14.5,2.3)
 \qbezier(12.5,1.8)(12.5,2.6)(12.5,2.6)
 \qbezier(12.8,2.9)(12.5,2.9)(14.5,2.9)
 \qbezier(12.5,2.6)(12.5,2.9)(12.8,2.9)
 \qbezier(14.5,2.9)(15.5,2.9)(14.5,2.3)


 \qbezier(.5,3.6)(.5,2.6)(2.5,4.1)
 \qbezier(.5,3.6)(.5,4.4)(.5,4.4)
 \qbezier(.8,4.7)(.5,4.7)(2.5,4.7)
 \qbezier(.5,4.4)(.5,4.7)(.8,4.7)
 \qbezier(2.5,4.7)(3.5,4.7)(2.5,4.1)

 \qbezier(3.5,3.6)(3.5,2.6)(5.5,4.1)
 \qbezier(3.5,3.6)(3.5,4.4)(3.5,4.4)
 \qbezier(3.8,4.7)(3.5,4.7)(5.5,4.7)
 \qbezier(3.5,4.4)(3.5,4.7)(3.8,4.7)
 \qbezier(5.5,4.7)(6.5,4.7)(5.5,4.1)

 \qbezier(6.5,3.6)(6.5,2.6)(8.5,4.1)
 \qbezier(6.5,3.6)(6.5,4.4)(6.5,4.4)
 \qbezier(6.8,4.7)(6.5,4.7)(8.5,4.7)
 \qbezier(6.5,4.4)(6.5,4.7)(6.8,4.7)
 \qbezier(8.5,4.7)(9.5,4.7)(8.5,4.1)

 \qbezier(9.5,3.6)(9.5,2.6)(11.5,4.1)
 \qbezier(9.5,3.6)(9.5,4.4)(9.5,4.4)
 \qbezier(9.8,4.7)(9.5,4.7)(11.5,4.7)
 \qbezier(9.5,4.4)(9.5,4.7)(9.8,4.7)
 \qbezier(11.5,4.7)(12.5,4.7)(11.5,4.1)

 \qbezier(12.5,3.6)(12.5,2.6)(14.5,4.1)
 \qbezier(12.5,3.6)(12.5,4.4)(12.5,4.4)
 \qbezier(12.8,4.7)(12.5,4.7)(14.5,4.7)
 \qbezier(12.5,4.4)(12.5,4.7)(12.8,4.7)
 \qbezier(14.5,4.7)(15.5,4.7)(14.5,4.1)

 \qbezier(.5,5.5)(.5,4.5)(2.5,6)
 \qbezier(.5,5.5)(.5,6.3)(.5,6.3)
 \qbezier(.8,6.6)(.5,6.6)(2.5,6.6)
 \qbezier(.5,6.3)(.5,6.6)(.8,6.6)
 \qbezier(2.5,6.6)(3.5,6.6)(2.5,6)

 \qbezier(3.5,5.5)(3.5,4.5)(5.5,6)
 \qbezier(3.5,5.5)(3.5,6.3)(3.5,6.3)
 \qbezier(3.8,6.6)(3.5,6.6)(5.5,6.6)
 \qbezier(3.5,6.3)(3.5,6.6)(3.8,6.6)
 \qbezier(5.5,6.6)(6.5,6.6)(5.5,6)

 \qbezier(6.5,5.5)(6.5,4.5)(8.5,6)
 \qbezier(6.5,5.5)(6.5,6.3)(6.5,6.3)
 \qbezier(6.8,6.6)(6.5,6.6)(8.5,6.6)
 \qbezier(6.5,6.3)(6.5,6.6)(6.8,6.6)
 \qbezier(8.5,6.6)(9.5,6.6)(8.5,6)

 \qbezier(9.5,5.5)(9.5,4.5)(11.5,6)
 \qbezier(9.5,5.5)(9.5,6.3)(9.5,6.3)
 \qbezier(9.8,6.6)(9.5,6.6)(11.5,6.6)
 \qbezier(9.5,6.3)(9.5,6.6)(9.8,6.6)
 \qbezier(11.5,6.6)(12.5,6.6)(11.5,6)

 \qbezier(12.5,5.5)(12.5,4.5)(14.5,6)
 \qbezier(12.5,5.5)(12.5,6.3)(12.5,6.3)
 \qbezier(12.8,6.6)(12.5,6.6)(14.5,6.6)
 \qbezier(12.5,6.3)(12.5,6.6)(12.8,6.6)
 \qbezier(14.5,6.6)(15.5,6.6)(14.5,6)


 \qbezier(.5,7.3)(.5,6.3)(2.5,7.8)
 \qbezier(.5,7.3)(.5,8.1)(.5,8.1)
 \qbezier(.8,8.4)(.5,8.4)(2.5,8.4)
 \qbezier(.5,8.1)(.5,8.4)(.8,8.4)
 \qbezier(2.5,8.4)(3.5,8.4)(2.5,7.8)

 \qbezier(3.5,7.3)(3.5,6.3)(5.5,7.8)
 \qbezier(3.5,7.3)(3.5,8.1)(3.5,8.1)
 \qbezier(3.8,8.4)(3.5,8.4)(5.5,8.4)
 \qbezier(3.5,8.1)(3.5,8.4)(3.8,8.4)
 \qbezier(5.5,8.4)(6.5,8.4)(5.5,7.8)

 \qbezier(6.5,7.3)(6.5,6.3)(8.5,7.8)
 \qbezier(6.5,7.3)(6.5,8.1)(6.5,8.1)
 \qbezier(6.8,8.4)(6.5,8.4)(8.5,8.4)
 \qbezier(6.5,8.1)(6.5,8.4)(6.8,8.4)
 \qbezier(8.5,8.4)(9.5,8.4)(8.5,7.8)

 \qbezier(9.5,7.3)(9.5,6.3)(11.5,7.8)
 \qbezier(9.5,7.3)(9.5,8.1)(9.5,8.1)
 \qbezier(9.8,8.4)(9.5,8.4)(11.5,8.4)
 \qbezier(9.5,8.1)(9.5,8.4)(9.8,8.4)
 \qbezier(11.5,8.4)(12.5,8.4)(11.5,7.8)

 \qbezier(12.5,7.3)(12.5,6.3)(14.5,7.8)
 \qbezier(12.5,7.3)(12.5,8.1)(12.5,8.1)
 \qbezier(12.8,8.4)(12.5,8.4)(14.5,8.4)
 \qbezier(12.5,8.1)(12.5,8.4)(12.8,8.4)
 \qbezier(14.5,8.4)(15.5,8.4)(14.5,7.8)

 \qbezier(.5,9.1)(.5,8.1)(2.5,9.6)
 \qbezier(.5,9.1)(.5,9.9)(.5,9.9)
 \qbezier(.8,10.2)(.5,10.2)(2.5,10.2)
 \qbezier(.5,9.9)(.5,10.2)(.8,10.2)
 \qbezier(2.5,10.2)(3.5,10.2)(2.5,9.6)

 \qbezier(3.5,9.1)(3.5,8.1)(5.5,9.6)
 \qbezier(3.5,9.1)(3.5,9.9)(3.5,9.9)
 \qbezier(3.8,10.2)(3.5,10.2)(5.5,10.2)
 \qbezier(3.5,9.9)(3.5,10.2)(3.8,10.2)
 \qbezier(5.5,10.2)(6.5,10.2)(5.5,9.6)

 \qbezier(6.5,9.1)(6.5,8.1)(8.5,9.6)
 \qbezier(6.5,9.1)(6.5,9.9)(6.5,9.9)
 \qbezier(6.8,10.2)(6.5,10.2)(8.5,10.2)
 \qbezier(6.5,9.9)(6.5,10.2)(6.8,10.2)
 \qbezier(8.5,10.2)(9.5,10.2)(8.5,9.6)

 \qbezier(9.5,9.1)(9.5,8.1)(11.5,9.6)
 \qbezier(9.5,9.1)(9.5,9.9)(9.5,9.9)
 \qbezier(9.8,10.2)(9.5,10.2)(11.5,10.2)
 \qbezier(9.5,9.9)(9.5,10.2)(9.8,10.2)
 \qbezier(11.5,10.2)(12.5,10.2)(11.5,9.6)

 \qbezier(12.5,9.1)(12.5,8.1)(14.5,9.6)
 \qbezier(12.5,9.1)(12.5,9.9)(12.5,9.9)
 \qbezier(12.8,10.2)(12.5,10.2)(14.5,10.2)
 \qbezier(12.5,9.9)(12.5,10.2)(12.8,10.2)
 \qbezier(14.5,10.2)(15.5,10.2)(14.5,9.6)


 \put(23,6){\vector(-1,0){7.5}}
 \qbezier(15,5.5)(19,3.2)(23,5)
 \put(22.5,4.3){\vector(1,1){1}}

 \put(18.5,6.6){$x=\epsilon\; y$}
 \put(18.5,3.2){$y=\frac{x}{\epsilon}$}


 \put(23,4){\vector(1,0){10}}
 \put(24,2){\vector(0,1){8}}


 \put(33,3){$y_1$}
 \put(22.5,9.5){$y_2$}

 \put(10,12){$D_\epsilon=D\backslash S_\epsilon$}
 \put(10,12){\vector(-2,-1){3}}

 \put(29,2.8){$l_1$}
 \put(23,6.5){$l_2$}
 \put(25,5){$S$}
 \put(30,5.5){\vector(-2,-1){1.7}}
 \put(30,5){$Y^*=Y\backslash\overline{S}$}
 \put(24.5,7.5){$Y=[0,l_1)\times[0,l_2)$}
 \put(22.8,2.8){$O$}

 \qbezier(24,7)(29,7)(29,7)
 \qbezier(29,4)(29,7)(29,7)

 \qbezier(24.5,5.3)(24.5,3.3)(28,6)
 \qbezier(24.5,5.3)(24.5,6.3)(24.5,6.3)
 \qbezier(24.8,6.6)(24.5,6.6)(28,6.6)
 \qbezier(24.5,6.3)(24.5,6.6)(24.8,6.6)
 \qbezier(28,6.6)(29,6.6)(28,6)

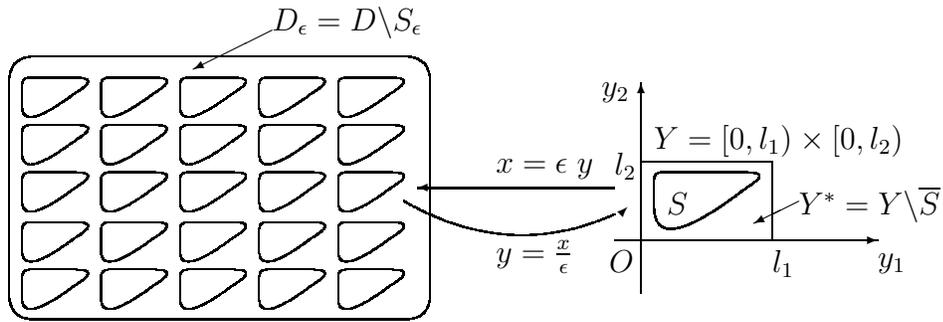
\put(10,-2.5){Fig. 1: Geometric setup in $\mathbb{R}^2$}

\end{picture}
\end{center}

\vskip 1.2cm

Now for $T>0$ fixed final time,  we consider the following
It$\hat{o}$ type nonautonomous stochastic partial differential
equation defined on the perforated domain $D_\e$ in
$\mathbb{R}^n$.
\begin{eqnarray}\label{e1}
du_\e(x,t)&=&\Big(div\big(A_\e(x)\nabla u_\e(x,t)\big)+
f_\e(x,t)\Big)dt+g_\e(t) dW(t)\\&& in
\;\; D_\e\times (0, T),\nonumber\\
u_\e&=&0\;\; on\;\; \p D\times (0, T),\\
\frac{\p u_\e}{\p \nu_{A_\e}}&=&0\;\; on \;\; \p S_\e\times (0,
T),\\
u_\e(0)&=&u_\e^0\;\; in \;\;  D_\e \label{e4},
\end{eqnarray}
where the matrix $A_\e$ is
\begin{equation*}
A_\e=\Big( a_{ij}\Big(\frac{x}{\e}\Big)\Big)_{ij}
\end{equation*}
and
\begin{equation*}
\frac{\p \;\cdot}{\p \nu_{A_\e}}=\sum_{ij}a_{ij}\Big(\frac{x}{\e}
\Big)\frac{\p\;\cdot}{\p x_j}n_i
\end{equation*}
with $n$ the exterior unit normal vector on the boundary $\p
D_\e$.

We make the following assumptions on the coefficients:
\begin{enumerate}
    \item $a_{ij}\in L^\infty(\R^n)$,\;\;$i,j=1,\cdots,n;$
    \item  $\sum_{i,j=1}^n a_{ij}\xi_i\xi_j \geq \alpha \sum_{i=1}^n\xi_i^2$ for
    $\xi\in\R^n$ and $\alpha$ a positive constant;
    \item $a_{ij}$ are $Y$-periodic.
\end{enumerate}
Furthermore we assume that

\begin{equation}\label{f}
f_\e\in L^2(D_\e\times [0, T])
\end{equation}
 and for $0\leq t\leq T$,
 $g_\e(t)$ is  a linear operator from $\ell^2$ to $
 L^2(D_\e)$ defined as
\begin{equation*}
g_\e(t) k=\sum_{i=1}^\infty g_\e^i(x,t) k_i,\;\; k=(k_1, k_2,\cdots
)\in\ell^2
\end{equation*}
where $g_\e^i(x,t)\in L^2(D_\e\times [0, T])$, $i=1,2,\cdots$, are
measurable functions with
\begin{equation}\label{ge}
\sum_{i=1}^\infty |g_\e^i(x,t)|^2_{L^2(D_\e)}< C_T, \;\;t\in [0, T]
\end{equation}
for some positive constant $C_T$ independent of $\e$. In (\ref{e1}),
$W(t)=(W_1(t), W_2(t), \cdots )$ is a Wiener process in $\ell^2$
with covariance operator $Q=Id_{\ell^2}$ and $\{W_i(t):
i=1,2,\cdots\}$ are mutually independent real valued standard Wiener
processes on a complete probability space $(\Omega, \mathcal{F},
\mathbb{P})$ with a canonical filtration $(\mathcal{F}_t)_{t\geq
0}$. Then
\begin{equation}\label{g}
|g_\e(t)|^2_{\mathcal{L}_2^Q}=\sum_{i=1}^\infty
|g_\e^i(x,t)|^2_{L^2(D_\e)}<C_T,\;\;t\in [0,T].
\end{equation}
Here $\mathcal{L}_2^Q$ is the space of Hilbert-Schmit operators
\cite{PZ92, HuangYan}. Denote by $\mathbf{E}$ the expectation
operator with respect to $\mathbb{P}$.

The following compactness result \cite{Li69} will be used in our
approach. Let $\mathcal{X}\subset \mathcal{Y}\subset \mathcal{Z}$
be three reflective Banach spaces and $\mathcal{X}\subset
\mathcal{Y}$ with compact and dense embedding. Define Banach space
$$
G=\{v: v\in L^2(0 ,T; \mathcal{X}), \frac{dv}{dt}\in L^2(0, T;
\mathcal{Z}) \}
$$
with norm
$$
|v|^2_G=\int_0^T|v(s)|^2_\mathcal{X}ds+\int_0^T\Big|\frac{dv}{ds}(s)\Big|^2_\mathcal{Z}ds,\;\;
v\in G.
$$

\begin{lemma}\label{cmpt}
If $B$ is bounded in $G$, then it is precompact in $L^2(0, T;
\mathcal{Y})$.
\end{lemma}

\medskip

Let $\mathcal{S}$ be a Banach space and $\mathcal{S'}$ be the
strong dual space of $\mathcal{S}$. We recall the definitions and
some properties of weak convergence and $\textrm{weak}^*$
convergence \cite{Yosida}.
\begin{defn}\label{weak1}
A sequence $\{s_n\}$ in $\mathcal{S}$ is said to converge weakly
to $s\in\mathcal{S}$ if $\forall s'\in\mathcal{S'} $,
$$
\lim_{n\rightarrow\infty}(s',s_n)_{\mathcal{S'}, \mathcal{S}}=(s',
s)_{\mathcal{S'}, \mathcal{S}}
$$
which is written as $s_n\rightharpoonup s$ weakly in
$\mathcal{S}$. Note that $(s',s)$ denotes the value of the
continuous linear functional $s'$ at the point $s$.
\end{defn}

\begin{lemma}(\textbf{Eberlein-Shmulyan})
Assume that $\mathcal{S}$ is reflexive and let $\{s_n\}$ be a
bounded sequence in $\mathcal{S}$. Then there exists a subsequence
$\{s_{n_k}\}$ and $s\in\mathcal{S}$ such that
$s_{n_k}\rightharpoonup s$ weakly in $\mathcal{S}$ as
$k\rightarrow\infty$. If all the weak convergent subsequence of
$\{s_n\}$ has the same limit $s$, then the whole sequence
$\{s_n\}$ weakly converges to $s$.
\end{lemma}

\begin{defn}\label{weak*}
A sequence $\{s'_n\}$ in $\mathcal{S'}$ is said to converge
$weakly^*$ to $s'\in\mathcal{S'}$ if $\forall s\in\mathcal{S} $,
$$
\lim_{n\rightarrow\infty}(s'_n, s)_{\mathcal{S'},
\mathcal{S}}=(s', s)_{\mathcal{S'}, \mathcal{S}}
$$
which is written as $s'_n\rightharpoonup s'$ $weakly^*$ in
$\mathcal{S'}$.
\end{defn}

\begin{lemma}
Assume that the dual space $\mathcal{S'}$ is reflexive and let
$\{s'_n\}$ be a bounded sequence in $\mathcal{S'}$. Then there
exists a subsequence $\{s'_{n_k}\}$ and $s'\in\mathcal{S'}$ such
that $s'_{n_k}\rightharpoonup s'$ $weakly^*$ in $\mathcal{S'}$ as
$k\rightarrow\infty$. If all the $weakly^*$ convergent subsequence
of $\{s_n'\}$ has the same limit $s'$, then the whole sequence
$\{s'_n\}$ $waekly^*$ converges to $s'$.
\end{lemma}

We also use the following definition of the weak convergence of
the Borel probability measures on $\mathcal{S}$, for more we refer
to \cite{Dudley}.

\begin{defn}\label{weak2}
Let $\{\mu_\e\}_\e$ be a family of Borel probability measures on
the Banach space $\mathcal{S}$. We say $\mu_\e$ weakly converges
to a Borel measure $\mu$ on $\mathcal{S}$ if
$$
\int_{\mathcal{S}}hd\mu_\e\rightarrow \int_{\mathcal{S}}hd\mu,\;\;
as\;\;\e\downarrow 0,
$$
for any $h\in C_b(\mathcal{S})$, the space of bounded continuous
functions on $\mathcal{S}$.
\end{defn}

In the following, for a fixed $T>0$, we always denote by $C_T$ a
constant
independent of $\e$.\\


\section{Basic properties of  the microscopic model}\label{s3}
In this section we will present some estimates of the solutions of
microscopic model (\ref{e1}), useful for the tightness result of the
distributions of solution processes in some appropriate space.

Let $H=L^2(D)$ and $H_\e=L^2(D_\e)$. Define the following space
\begin{equation*}
V_\e=\{u\in H^1(D_\e), u|_{\p D}=0 \}
\end{equation*}
provided with the norm
$$
|v|_{V_\e}=|\nabla_{A_\e}v |_{\oplus_n H_\e}=\Big
|\Big(\sum_{j=1}^na_{ij}\Big(\frac{x}{\e} \Big)\frac{\p v }{\p
x_j}\Big)_{i=1}^n \Big |_{\oplus_nH_\e}.
$$
This norm is equivalent to the usual $H^1(D_\e)$-norm, with an
embedding constant independent of $\e$,  due to the assumptions on
$a_{ij}$ in the last section. Here $\oplus_n$ denotes the direct
sum of the Hilbert spaces with usual direct sum norm.
 Let
$$\mathcal{D}(\mathcal{A}_\e)=\Big\{v\in V_\e: div(A_\e \nabla v)\in H_\e \;\;and
\;\; \frac{\p v}{\p \nu_{A_\e}}\Big|_{\p S_\e}=0 \Big\}
$$
and define operator $\mathcal{A}_\e v=div(A_\e\nabla v)$ for $v\in
\mathcal{D}(\mathcal{A}_\e)$. Then system (\ref{e1})-(\ref{e4})
can be written as the following abstract stochastic evolutionary
equation
\begin{equation}\label{abse}
du_\e=(\mathcal{A}_\e u_\e+f_\e)dt+g_\e dW,\;\; u_\e(0)=u_\e^0.
\end{equation}

By the assumptions on $a_{ij}$, operator $\mathcal{A}_\e$ generates
a  strongly continuous semigroup $S_\e(t)$ on $H_\e$.
 Solution of (\ref{abse}) can then be written in the mild sense
\begin{equation}\label{mild}
u_\e(t)=S_\e(t)u_\e^0+\int_0^tS_\e(t-s)f_\e(s)ds+\int_0^tS_\e(t-s)g_\e(s)dW(s)
\end{equation}
And the variational formulation is
\begin{eqnarray}\label{variational}
\big(d u_\e(t),
v\big)_{H_\e^{-1},V_\e}&=&\Big(-\int_{D_\e}A_\e(x)\nabla
u_\e(x,t)\nabla v(x)dx+\int_{D_\e}f_\e(x,t)v(x)dx\Big)dt+\nonumber \\
&&\int_{D_\e}g_\e(x,t)v(x)dW(t),\;\; in \;\; \mathcal{D}'(0,
T),\;\; v\in V_\e,
\end{eqnarray}
with $u_\e(0,x)=u^0_\e(x)$. \\

For the well-posedness of  system (\ref{abse}) we have the following
result.


\begin{theorem}\label{wellpose} (\textbf{Global well-posedness of microscopic model})
Assume (\ref{f}) and (\ref{g}) hold. Let $u_\e^0$ be a
$\big(\mathcal{F}_0, \mathcal{B}(H_\e)\big)$-measurable random
variable. Then system (\ref{abse}) has a unique mild solution
$u\in L^2\big(\Omega, C(0,T; H_\e)\cap L^2(0, T; V_\e)\big)$,
which is also a weak solution in the following sense
\begin{eqnarray}
&&\hspace{-0.3cm} (u_\e(t),v)_{H_\e}\nonumber\\
\hspace{-0.3cm}&=&\hspace{-0.3cm}(u_\e^0, v)_{H_\e}+\int_0^t(\mathcal{A}_\e u_\e(s),
v)_{H_\e}ds+\int_0^t(f_\e, v)_{H_\e}ds+ \int_0^t(g_\e dW,
v)_{H_\e}\label{weak}
\end{eqnarray}
for $t\in[0,T)$ and $v\in V_\e$. Moreover, if $u_\e^0$ is
independent of   $W(t)$ with $\mathbf{E}|u_\e^0|^2_{H_\e} <\infty
$, then
\begin{equation}\label{est1}
\mathbf{E}|u_\e(t)|^2_{H_\e}+\mathbf{E}\int_0^t|u_\e(s)|^2_{V_\e}ds\leq
\mathbf{E}|u_\e^0|^2_{H_\e}+C_T,\;\; for \;\; t\in [0, T],
\end{equation}
\begin{equation}\label{est2}
\mathbf{E}\int_0^t|\dot{u}_\e(s)|^2_{H^{-1}_\e}ds\leq
C_T(\mathbf{E}|u_\e^0|^2_{H_\e}+1),\;\; for \;\; t\in [0, T].
\end{equation}
If further assume that
\begin{equation}\label{g1}
|\nabla_{A_\e}g_\e(t)|^2_{\mathcal{L}_2^Q}=\sum_{i=1}^\infty|\nabla_{A_\e}g_\e^i(t)|^2_{\oplus_nH_\e}\leq
C_T,\;\;for \;t\in [0, T]
\end{equation}
with $u_\e^0\in V_\e$ and $\mathbf{E}|u_\e^0|^2_{V_\e}<\infty$,
then

\begin{equation}\label{est3}
\mathbf{E}|u_\e(t)|^2_{V_\e}+\mathbf{E}\int_0^t|\mathcal{A}_\e
u_\e(s)|^2_{H_\e}ds\leq \mathbf{E}|u_\e^0|^2_{V_\e}+C_T,\;\; for
\;\; t\in [0, T].
\end{equation}

Moreover, system (\ref{abse}) is well-posed on $[0, \infty)$ if
\begin{equation}\label{fg}
f_\e\in L^2(0,\infty; H_\e), \;\; g_\e\in L^2(0, \infty;
\mathcal{L}_2^Q).
\end{equation}
\end{theorem}

\begin{proof}
By the assumption (\ref{g}), we have
$$
\|g_\e(t)\|^2_{\mathcal{L}^Q_2}=\sum^\infty_{i=1}|g^i_\e(t,x)|^2_{H_\e}<\infty.
$$
Then the classical result of \cite{PZ92} yields the local existence
of $u_\e$. And applying the stochastic Fubini theorem, it is easy to
verify the local mild solution is also a weak solution.

Now we give the following a  priori  estimates which yields the
existence of weak solution on $[0 ,T]$ provide (\ref{f}) and
(\ref{g}) hold.

Applying It$\hat{o}$ formula to $|u_\e|^2$, we derive
\begin{equation}\label{It1}
d|u_\e(t)|^2_{H_\e}-2(\mathcal{A}_\e u_\e, u_\e)_{H_\e}dt=2(f_\e,
u_\e)_{H_\e}dt+2(g_\e dW, u_\e)_{H_\e}+|g_\e|^2_{\mathcal{L}_2^Q}dt.
\end{equation}
By the assumption on $a_{ij}$, we see that
$$
-(\mathcal{A}_\e u_\e, u_\e)_{H_\e}\geq\lambda |u_\e|^2_{H_\e}
$$
for some constant $\lambda>0$ independent of $\e$. Then
 integrating (\ref{It1}) with respect to $t$ yields
\begin{eqnarray*}
&&|u_\e(t)|^2_{H_\e} +\int_0^t|u_\e|^2_{V_\e}ds\\ &\leq
&|u_\e^0|^2_{H_\e}+\lambda^{-1}|f_\e|^2_{L^2(0,T;
H_\e)}+\int_0^t(g_\e dW,
u_\e)_{H_\e}ds+\int_0^t|g_\e|^2_{\mathcal{L}_2^Q}ds.
\end{eqnarray*}
Taking expectation on both sides of the above inequality, we
derive (\ref{est1}).

In a similar way, application of  It$\hat{o}$ formula to
$|u_\e|^2_{V_\e}=|\nabla_{A_\e}u_\e|^2_{\oplus_n H_\e}$ results in the
relation
\begin{eqnarray}
&&d|u_\e(t)|^2_{V_\e} +2(\mathcal{A}_\e u_\e, \mathcal{A}_\e
u_\e)_{H_\e}dt \nonumber\\&=&-2(f_\e,\mathcal{A}_\e
u_\e)_{H_\e}dt-2(g_\e dW, \mathcal{A}_\e u_\e)_{H_\e}+
|\nabla_{A_\e}g_\e|^2_{\mathcal{L}_2^Q}dt. \label{It2}
\end{eqnarray}
Integrating both sides of (\ref{It2}) and by the Cauchy-Schwarz
inequality, it is easily to have
\begin{eqnarray*}
&&|u_\e(t)|_{V_\e}^2+\int_0^t|\mathcal{A}_\e u_\e|^2_{H_\e}ds
\\ &\leq& |u_\e(0)|^2_{V_\e}+|f_\e|^2_{L^2(0,T; H_\e)}-2\int_0^t (g_\e dW,
\mathcal{A}_\e u_\e)_{H_\e}ds+\int_0^t
|\nabla_{A_\e}g_\e|^2_{\mathcal{L}_2^Q}ds,
\end{eqnarray*}
Then taking the expectation, we derive (\ref{est3}). By
(\ref{variational}) and the property of the stochastic integral we
easily have (\ref{est2}).

Thus, by the above estimates, the solution can be extended to
$[0,\infty)$ if (\ref{fg}) hold. The proof is complete.
\end{proof}

\medskip

We recall a probability concept. Let $z$ be a random variable
taking values in a Banach space $\mathcal{S}$, namely, $z:\; \Om
\to z$. Denote by $\mathcal{L}(z)$ the distribution (or law) of
$z$. In fact, $\mathcal{L}(z)$ is a Borel probability measure on
$\mathcal{S}$ defined as \cite{PZ92}

$$
\mathcal{L}(z)(A)=\mathbb{P} \{\omega: z(\omega) \in A \},
$$
for every event (i.e., a Borel set) $A$ in the Borel
$\sigma-$algebra $\mathcal{B}(\mathcal{S})$, which is the smallest
$\sigma-$algebra containing all open balls in $\mathcal{S}$.

\vskip 0.5cm

 As stated in \S   \ref{s1},  for the SPDE (\ref{e1}) we aim at deriving
 an
 effective equation in the sense of probability. A solution  $u_\e$ may be regarded as a
   random variable  taking values in  $L^2(0, T; H_\e)$.
  So for a solution $u_\e$ of
 (\ref{e1})-(\ref{e4}) defined on $[0, T]$,  we
focus on the behavior of distribution of $u_\e$ in $L^2(0, T;
H_\e)$ as $\e\rightarrow 0$. For this purpose, the tightness
\cite{Dudley} of
  distributions is needed. Note   that the
function space changes with $\e$, which is a difficulty for
obtaining  the tightness of   distributions.   Thus we will treat
$\{\mathcal{L}(u_\e)\}_{\e>0}$ as a family of distributions on
$L^2(0, T; H)$ by extending $u_\e$ to the whole domain $D$. Recall
that the distribution (or law ) of  $u_\e$ is defined as:
$$
\mathcal{L}(u_\e)(A)=\mathbb{P} \{\omega: u_\e(\cdot, \cdot, \omega)
\in A \}
$$
for Borel set $A$ in $ L^2(0, T;H_\e)$. First we define an extension
operator $P_\e$ in the following lemmas.

\medskip

In the following we denote by $\mathbb{L}\big(\mathcal{X},
\mathcal{Y})$ the space of bounded linear  operator from Banach
space $\mathcal{X}$ to Banach space $\mathcal{Y}$.

\begin{lemma}\label{Q}
There exists a bounded linear operator
$$
 \hat{Q}\in\mathbb{L}(H^k(Y^*), H^k(Y)),\;\;k=0, 1,
$$
such that
$$
|\nabla \hat{Q} v|_{\oplus_n L^2(Y)}\leq C |\nabla v|_{\oplus_n
L^2(Y^*)},\;\; v\in H^1(Y^*)
$$
for some constant $C>0$.
\end{lemma}
For the proof of Lemma \ref{Q} see \cite{CD89}.

 We define an extension operator $P_\e$ in terms of the above bounded
 linear operator $\hat{Q}$ in the following lemma.
\begin{lemma}\label{Pe}
There exists an extension operator
$$
 P_\e\in \mathbb{L}\big(L^2(0,T; H^k(D_\e)), L^2(0,T;
H^k(D))\big),\;\;k=0,1,
$$
such that for any $v\in H^k(D_\e)$
\begin{enumerate}
    \item $P_\e v=v $\;\;on $D_\e\times (0,T)$
    \item $|P_\e v|_{L^2(0,T; H)}\leq C_T|v|_{L^2(0,T; H_\e)}$
    \item $|\nabla_{A_\e}(P_\e v)|_{L^2(0, T; \oplus_n L^2(D))}\leq C_T|\nabla_{A_\e} v|_{L^2(0,T;\oplus_n L^2(D_\e) )}$
\end{enumerate}
where $C_T$ is a constant independent of $\e$.
\end{lemma}
\begin{proof}
For $\phi\in H^k(D_\e)$, then
$$
\phi_\e(y)=\frac{1}{\e}\phi(\e\;y)
$$
belongs to $H^k(Y^*_l)$ with $Y^*_l$ the translation of $Y^*$ for
some $l\in \R^n$. Define
\begin{equation}\label{Qe}
\hat{Q}_\e\phi(x)=\e(\hat{Q}\phi_\e)\big(\frac{x}{\e}\big).
\end{equation}
Now for $\phi\in L^2(0, T;H^k(D_\e))$, we define
$$
(P_\e\phi)(x,t)=[\hat{Q}_\e\phi(t,\cdot)]\Big(\frac{x}{\e}\Big)
=\e[\hat{Q}\phi_\e(t,\cdot)]\Big(\frac{x}{\e}\Big).
$$
It is known  \cite{CD89}  that the operator $P_\e\in
\mathbb{L}\big(L^2(0,T; H^k(D_\e)), L^2(0,T; H^k(D))\big),\;\;k$
\\ $=0,1,$ and satisfies the conditions (1)-(3) listed in the lemma.
This completes the proof.
\end{proof}

\begin{remark}
In Lemma 2.1 of \cite{CD89}, the operator $P_\e$ defined in
 $\mathbb{L}\big(L^\infty(0,T;$\\ $ H^k(D_\e)),$  $ L^\infty(0,T;
H^k(D))\big),\;\;k=0,1,$ coincides with the operator defined in
Lemma \ref{Pe} above.
\end{remark}

\begin{remark}\label{rem1}
The  estimates  in Theorem  \ref{wellpose} for $u_\e$ also hold
for $P_\e u_\e$. In fact estimates (\ref{est1}) and (\ref{est3})
are easily derived due to the property of the operator of $P_\e$.
Since the operator $P_\e$ is defined on $L^2(0,T; H^k(D_\e)$,
$k=0,1$, for (\ref{est2}) we define
$$
P_\e\dot{u}_\e\equiv\mathcal{A}_\e P_\e
u_\e+\tilde{f}_\e+\tilde{g}_\e\dot{W},\;\;on\; D\times (0, T).
$$
By the property of $P_\e$ and the estimates of $u_\e$, it is easy
to see that
$$
P_\e\dot{u}_\e=\dot{(P_\e u_\e)},\;\;in \;D_\e\times (0, T)
$$
and
$$
\mathbf{E}|P_\e \dot{u}_\e|_{L^2(0,T;H^{-1})}\leq
\mathbf{E}|\dot{u}_\e|_{L^2(0,T;H^{-1}_\e)}.
$$
\end{remark}


\section{Effective macroscopic model}\label{s4}

We now derive the effective macroscopic model for the original
model (\ref{e1}). Let $u_\e\in L^2(0, T; H_\e )$ be the solution
of system (\ref{e1})-(\ref{e4}). Then by the estimates in Theorem
\ref{wellpose}, Remark \ref{rem1} and the Chebyshev inequality
\cite{PZ92, Dudley}, it is clear that for any $\delta>0$ there is
a bounded set $K_\delta\subset G$ with spaces $\mathcal{X}$,
$\mathcal{Y}$ and $\mathcal{Z}$ in Lemma \ref{cmpt} (and in the
paragraph immediately before it) are replaced by $H_0^1(D)$, $H$
and $H^{-1}(D)$ respectively, such that
\begin{equation*}
\mathbb{P}\{P_\e u_\e\in K_\delta \}>1-\delta.
\end{equation*}
Thus $K_\delta$ is compact in $L^2(0,T; H)$ by Lemma \ref{cmpt}.
Then $\{\mathcal{L}(P_\e u_\e)\}_\e$ is tight in $L^2(0, T; H)$.
The Prokhorov Theorem and the Skorohod embedding theorem
(\cite{PZ92}) assure  that for any sequence $\{\e_j\}$ with
$\e_j\rightarrow 0$ as $j\rightarrow \infty$, there exists a
subsequence $\{\e_{j(k)}\}$, random variables
$\{\hat{u}_{\e_{j(k)}}\}\subset L^2(0, T; H_{\e_{j(k)}})$ and
$u\in L^2(0, T; H)$  defined on a new probability space
$(\widehat{\Omega}, \widehat{\mathcal{F}}, \widehat{\mathbb{P}})$,
such that
$$
\mathcal{L}(P_{\e_{j(k)}}\hat{u}_{\e_{j(k)}})=\mathcal{L}(P_{\e_{j(k)}}u_{\e_{j(k)}})
$$
and
$$
P_{\e_{j(k)}}\hat{u}_{\e_{j(k)}}\rightarrow u\;\;in\;\;L^2(0,T;
H)\;\; as \;\;k\rightarrow \infty,
$$
for almost all $\omega\in\widehat{\Omega}$. Moreover
$P_{\e_{j(k)}}\hat{u}_{\e_{j(k)}}$ solves system
(\ref{e1})-(\ref{e4}) with $W$ replaced by Wiener process
$\widehat{W}_k$ defined on probability space $(\widehat{\Omega},
\widehat{\mathcal{F}}, \widehat{\mathbb{P}})$ with same
distribution as $W$. The limit $u$ is unique; see \cite{Bill2},
p.333. In the following, we will determine the limiting equation
(homogenized effective equation) that $u$ satisfies and  the
limiting equation is independent of $\e$. After this is done we
see that $\mathcal{L}(u_\e)$ weakly converges to $\mathcal{L}(u)$ as $\e\downarrow 0$.\\

 We always
assume the following
\begin{equation}\label{fe}
\tilde{f}_\e\rightharpoonup f,\;\;weakly\; in\;\;L^2(0,T;H)
\end{equation}
and
\begin{equation}\label{gee}
\tilde{g}^i_\e\rightharpoonup g^i,\;\; weakly\; in\;\;L^2(0,T;H).
\end{equation}
Define a new probability space $(\Omega_\delta,\mathcal{F}_\delta,
\mathbb{P}_\delta)$ as
\begin{equation*}
\Omega_\delta=\{\omega\in\Omega: u_\e(\omega)\in K_\delta\},
\end{equation*}
\begin{equation*}
\mathcal{F}_\delta=\{F\cap\Omega_\delta: F\in\mathcal{F}\}
\end{equation*}
and
$$
\mathbb{P}_\delta(F)=\frac{\mathbb{P}(F\cap\Omega_\delta)}{\mathbb{P}(\Omega_\delta)},
\;\;for \;\; F\in\mathcal{F}_\delta.
$$
Denote by $\mathbf{E}_\delta$ the expectation operator with respect
to $\mathbb{P}_\delta$.

Now we restrict the system on the  probability space
$(\Omega_\delta, \mathcal{F}_\delta, \mathbb{P}_\delta)$. In the
following discussion  we aim at obtaining $L^2(\Omega_\delta)$
convergence for any $\delta>0$ which means the convergence in
probability \cite{Bill, Dudley}.

From the estimates $(\ref{est1})$, $(\ref{est2})$, Remark
\ref{rem1} and the compact injection $G\subset L^2(0, T; H)$,
there exists a subsequence of $u_\e$ in $K_\delta$, still denoted
by $u_\e$, such that for a fixed $\omega\in\Omega_\delta$
\begin{eqnarray}
&&P_\e u_\e\rightharpoonup u \;\; weakly^* \; in \;\;L^\infty(0,T; H)\\
&&P_\e u_\e\rightharpoonup u \;\; weakly  \; in \;\;L^2(0,T; H^1)\\
&&P_\e u_\e\rightarrow u \;\; strongly \; in \;\;L^2(0,T; H)\\
&&P_\e \dot{u}_\e\rightharpoonup \dot{u} \;\; weakly \; in
\;\;L^2(0,T; H^{-1}).
\end{eqnarray}

Define
$$
\xi_\e=\Big(\sum^n_{j=1}a_{ij}\Big(\frac{x}{\e}\Big)\frac{\p
u_\e}{\p x_j} \Big)=A_\e\nabla u_\e
$$
which satisfies
\begin{eqnarray}
-div\xi_\e&=&f_\e+g_\e\dot{W}-\dot{u}_\e\;\; in\; D_\e\times(0, T) \\
\xi_\e\cdot n &=&0\;\;on\; \p S_\e\times(0, T).
\end{eqnarray}
By the hypothesis of $a_{ij}$ and the fact that
$(\tilde{u}_\e)_\e$ being bounded in $L^2(0, T; H^1_0)$, we have
\begin{equation}\label{conv1}
\tilde{\xi}_\e\rightharpoonup \xi\;\; weakly\; in \;\; L^2(0, T;
\oplus_n H).
\end{equation}
We make use of Tartar's method of oscillating test functions to
determine the limiting equation \cite{CD99}.

Note that
\begin{eqnarray}  \label{formula}
\int_0^T\int_D\tilde{\xi}_\e\cdot \nabla v \phi
dxdt&=&\int_0^T\int_D \tilde{f}_\e v\phi
dxdt+\sum_{i=1}^\infty\int_0^T\int_D\tilde{g}^i_\e v dx\phi dW_i(t)+\nonumber \\
&&\int_0^T\int_D P_\e u_\e\chi_{D_\e}\dot{\phi}vdxdt\label{xi}
\end{eqnarray}
for all $v\in H_0^1(D)$ and $\phi\in \mathcal{D}(0,T)$. We pass to
the limit in (\ref{xi}) as $\e\rightarrow 0$. Due to the fact
\begin{equation}\label{ue}
 P_\e u_\e \rightarrow u\;\;
\textrm{strongly}\;\;\textrm{in}\;\; L^2(0, T; H),
\end{equation}
\begin{equation}\label{chi}
\chi_{D_\e}\rightharpoonup \vartheta \;\; \textrm{weakly}^*\;\;
\textrm{in}\;\; L^\infty(D)
\end{equation}
and the estimate
$$
\mathbf{E}\Big|\sum_{i=1}^\infty\int_0^T\int_D\tilde{g}^i_\e v
dx\phi dW_i(t)\Big|^2\leq
\sum_{i=1}^\infty|\tilde{g}_\e^i|_{L^2(0, T;H)}^2|v\phi|^2_{L^2(0,
T; H)},
$$
by the assumption (\ref{gee}), we see that
$$
\sum_{i=1}^\infty\int_0^T\int_D\tilde{g}^i_\e v dx\phi
dW_i(t)\rightarrow \sum_{i=1}^\infty \int_0^T\int_D g^i v dx \phi
dW_i(t), \;\;in\;\; L^2(\Omega).
$$

Thus letting $\e\rightarrow 0$ in (\ref{xi}) and since $L^2(\Omega_\delta)$
is a subspace of $L^2(\Omega)$ one finds that in
$L^2(\Omega_\delta)$
\begin{eqnarray}\label{xi0}
\int_0^T\int_D\xi \cdot \nabla v \phi dxdt&=&\int_0^T\int_D f
v\phi dxdt+\sum_{i=1}^\infty \int_0^T\int_D g^i v dx \phi
dW_i(t)\nonumber \\ &&+\int_0^T\int_D \vartheta u \dot{\phi}vdxdt.
\end{eqnarray}
Hence
\begin{equation}\label{xi1}
-div \;
\xi(x,t)=f(x,t)+g(x,t)\dot{W}-\vartheta\dot{u}\;\;in\;\;D\times
(0, T).
\end{equation}

In the following we identify the limit $\xi$. We follow  the
approach of deterministic case for the elliptic problem  with
homogeneous Neumann boundary condition \cite{CD99}.

For any $\lambda\in\R^n$, let $w_\lambda$ be the solution of
\begin{eqnarray}
-\sum_{j=1}^n\frac{\p }{\p y_j}\Big ( \sum_{i=1}^n a_{ij}(y)\frac{\p
w_\lambda}{\p
y_i}\Big)=0 \;\; in \;\; Y^*\\
w_\lambda-\lambda\cdot y \;\;\;\; Y-periodic\\
\frac{\p w_\lambda}{\p \nu_A}=0\;\;\;\;\;\;\;\;\;\;\;\;\;\;\; on\;
\p S
\end{eqnarray}
and define
$$
w_\lambda^\e=\e(\hat{Q} w_\lambda)\Big(\frac{x}{\e}\Big)
$$
where $\hat{Q}$ is in Lemma \ref{Q}. Then we have \cite{CD99},
\begin{eqnarray}
w_\lambda^\e\rightharpoonup \lambda\cdot x\;\; weakly\;\;in\;\;H^1(D),\label{w0}\\
\nabla w_\lambda^\e \rightharpoonup \lambda\;\;
weakly\;\;in\;\;\oplus_n L^2(D).\label{w1}
\end{eqnarray}
Now we define
$$
(\eta^\lambda_j(y))_j=\Big(\sum_{i=1}^na_{ji}(y)\frac{\p
w_\lambda(y)}{\p y_i}\Big)_j,\;\;y\in Y^*
$$
and $(\eta^\lambda_\e)(x)=(\eta^\lambda_j(x/\e))_j= A^t_\e \nabla
w_\lambda^\e$. Then
\begin{equation}\label{eta0}
-div \; \tilde{\eta}^\lambda_\e=0\;\; in \;\; D
\end{equation}
and due to (\ref{w0}) and (\ref{w1})
\begin{equation}\label{eta1}
\tilde{\eta}^\lambda_\e\rightharpoonup
\mathcal{M}_Y(\eta^\lambda)\;\;weakly\;\;in \;\; L^2(D).
\end{equation}
It is easy to see that $\mathcal{M}_Y(\eta^\lambda)=B^t\lambda$
with $B^t=(\beta_{ji})$ a constant matrix which is determined in
the appendix.

Using test function  $\phi v w_\lambda^\e$ with $\phi\in
\mathcal{D}(0,T)$, $v\in\mathcal{D}(D)$ in (\ref{xi}) and
multiplying both sides of (\ref{eta0}) with  $\phi v P_\e u_\e$, one
has
\begin{multline*}
\int_0^T\int_D\tilde{\xi}_\e\cdot\nabla v \phi w_\lambda^\e
dxdt+\int_0^T\int_{D_\e}\xi_\e\cdot\nabla w_\lambda^\e v\phi dxdt \\
-\int_0^T\int_D\tilde{\eta}_\e^\lambda\cdot\nabla v\phi P_\e u_\e
dxdt-\int_0^T\int_D \tilde{\eta}_\e^\lambda\cdot\nabla(P_\e
u_\e)v\phi dxdt=
\end{multline*}
\begin{multline*}
 \int_0^T\int_D\tilde{f}_\e\phi v w_\lambda^\e
dxdt+\sum_{i=1}^\infty\int_0^T\int_D\tilde{g}^i_\e v w_\lambda^\e
dx\phi dW_i(t)+\int_0^T\int_D P_\e u_\e\chi_{D_\e}\dot{\phi} v
w_\lambda^\e dxdt.
\end{multline*}
Then by the definition of $\xi_\e$, $\eta_\e^\lambda$ and the
assumptions (\ref{fe}), (\ref{gee}), using the convergence
(\ref{conv1}), (\ref{ue}), (\ref{chi}), (\ref{w0}), (\ref{w1}) and
(\ref{eta1}), we have in $L^2(\Omega_\delta)$
\begin{multline*}
\int_0^T\int_D\xi\cdot\nabla v \phi \lambda\cdot x
dxdt-\int_0^T\int_DB^t\lambda\cdot\nabla v\phi u dxdt\\
=\int_0^T\int_Df\phi v \lambda\cdot
xdxdt+\sum_{i=1}^\infty\int_0^T\int_D g^i v \lambda\cdot x  dx\phi
dW_i(t)+\int_0^T\int_D\vartheta u v\dot{\phi}\lambda\cdot xdxdt.
\end{multline*}
That is
\begin{multline*}
\int_0^T\int_D\xi\cdot\nabla (v\lambda\cdot x) \phi dxdt-
\int_0^T\int_D\xi\cdot\lambda v\phi dxdt- \int_0^T\int_DB^t\lambda\cdot\nabla v\phi u dxdt\\
=\int_0^T\int_Df\phi v \lambda\cdot
xdxdt+\sum_{i=1}^\infty\int_0^T\int_D g^i v \lambda\cdot x  dx\phi
dW_i(t)+\int_0^T\int_D\vartheta u v\dot{\phi}\lambda\cdot xdxdt.
\end{multline*}
Then by using (\ref{xi0}) with the test function replaced by
$v\lambda\cdot x\phi$ one has
\begin{equation*}
\int_0^T\int_D\xi\cdot\lambda v\phi dxdt=
\int_0^T\int_DB^t\lambda\cdot\nabla u\phi v dxdt
\end{equation*}
which yields
\begin{equation*}
\xi\cdot\lambda=B^t\lambda\cdot\nabla u=B\nabla u\cdot \lambda.
\end{equation*}
Then
$$
\xi=B\nabla u
$$
since $\lambda$ is arbitrary. Then $u$ satisfies the following
equation
\begin{equation}\label{u}
\vartheta du=\big(div(B\nabla u)+f\big)dt+g dW(t).
\end{equation}

\vskip 0.3cm

Suppose
\begin{equation}\label{ue0}
\tilde{u}_\e^0\rightharpoonup u^0,\;\; weakly \;in\; H.
\end{equation}
We now determine the initial value by suitable test-functions.
 In fact, taking $v\in\mathcal{D}(D)$ and
$\phi\in\mathcal{D}([0,T])$ with $\phi(T)=0$ we have
\begin{eqnarray*}
\int_0^T\int_D\tilde{\xi}_\e\cdot\nabla v\phi
dxdt&=&\int_0^T\int_D\tilde{f}_\e v \phi
dxdt+\sum_{i=0}^\infty\int_0^T\int_D\tilde{g}_\e^ivdx\phi dW_i(t)
-\\ &&\int_0^T\int_D\tilde{u}_\e
v\dot{\phi}dxdt+\int_D\tilde{u}_\e^0 \phi(0)vdx.
\end{eqnarray*}
Now let $\e\rightarrow 0$, noticing that
\begin{eqnarray*}
\int_0^T\int_D\tilde{u}_\e
v\dot{\phi}dxdt&=&\int_0^T\int_D\chi_{D_\e}P_\e\tilde{u}_\e
v\dot{\phi}dxdt \rightarrow \int_0^T\int_D\vartheta u v
\dot{\phi}dxdt=\\&&-\int_0^T\int_D\vartheta \dot{u}v\phi
dxdt+\int_D\vartheta u(0)\phi(0)vdx
\end{eqnarray*}
by (\ref{xi1}), we have
$$
u(0)=\frac{u^0}{\vartheta}.
$$

\vskip 0.3cm

Here one should notice that the above result is in the sense of
$L^2(\Omega_\delta)$. Then  the above analysis yields  the
following results
\begin{equation}\label{r1}
\lim_{\e\rightarrow 0}\mathbf{E}_\delta |P_\e u_\e-u|^2_{L^2(0, T;
H)}=0
\end{equation}
and
\begin{equation}\label{r2}
\lim_{\e\rightarrow 0}\mathbf{E}_\delta
\int_0^T\int_D(\mathcal{A}_\e P_\e u_\e-B\nabla u)v\phi dxdt=0
\end{equation}
for any $v\in\mathcal{D}(D)$ and $\phi\in\mathcal{D}([0,T])$.\\

 Now we are in
the position to give the homogenized effective equation in the
following theorem.

\begin{theorem}\label{homo} (\textbf{Effective macroscopic model})
For any $T>0$, assume that (\ref{fe}), (\ref{gee}) and (\ref{ue0})
hold. Let $u_\e$ be the solution of (\ref{e1})-(\ref{e4}). Then
the distribution $\mathcal{L}(P_\e u_\e)$ converges weakly to
$\mu$ in the space of probability measures on $L^2(0, T; H)$ as
$\e\downarrow 0$, with $\mu$ being the distribution of $u$, which
is the solution of the following homogenized effective equation
\begin{eqnarray}
\vartheta du&=&\big(div(B\nabla u)+f\big)dt+gdW(t)\;\; in\;D \times (0, T), \label{limit} \\
u&=&0\;\; on \;\partial D\times (0,T),\\
u(x,0)&=&\frac{u^0}{\vartheta}\;\;in \;D,   \label{limit2}
\end{eqnarray}
where $B=(\beta_{ij})$ is determined by (\ref{beta}) in Appendix at
the end of this paper.
\end{theorem}

\begin{remark}
This theorem implies that the macroscopic  model (\ref{limit}) is
an effective approximation for the microscopic model (\ref{e1}),
on any finite time interval $0<t<T$, in the sense of probability
distribution. In other words, if we intend to numerically simulate
the  microscopic model up to finite time, we could use the
macroscopic model as an approximation when $\e$ is sufficiently
small.
\end{remark}

\begin{remark}
Due to the appearance of the stochastic integral term  (see
(\ref{formula})), this theorem on weak convergence of probability
measures does not follow directly from the deterministic
homogenization results and the mild formulation (\ref{mild}).
 \end{remark}

\begin{remark}\label{unnique}
The stochastic PDE   (\ref{limit}) is defined on the homogenized
domain $D$. By the   analysis  in \cite{PZ92}, for any fixed
$T>0$, the macroscopic system (\ref{limit})-(\ref{limit2}) is
well-posed, as long as  $f\in L^2(0, T; H)$ and $g\in L^2(0, T;
\mathcal{L}^Q_2)$.
\end{remark}
\medskip

\begin{proof}
Noticing the arbitrariness of $\delta$, this is direct result of the
analysis of the first part in this section by the Skorohod theorem
and the $L^2(\Omega_\delta)$ convergence of $P_\e u_\e$ on
$(\Omega_\delta, \mathcal{F}_\delta, \mathbb{P}_\delta)$.
\end{proof}

We finish this section by the following remark.
\begin{remark}
Note that there are several papers on   effective
 dynamics for      partial differential equations with random coefficients
 (so called random PDEs; not   stochastic PDEs); see
 \cite{KP02, PP03, Watanabe} and reference therein.
 In \cite{KP02, PP03}, a random partial
 differential equation is obtained as the homogenized effective equation
 for a random system with fast or small scales on both time or spatial variable.
 And the distribution of solution
 of heterogeneous system converges weakly to that of homogenized equation.
 However in \cite{Watanabe}, the effective  equation is obtained as an
averaged  deterministic equation for a random system with small
scale just on time. And the fluctuation of the solution of the
random equation around the solution of the averaged equation
converges to a generalized Ornstein-Uhlenbeck
 process in distribution. In the present paper, the original microscopic model
 is a stochastic PDE (i.e., PDE with white noise) and
 the     effective macroscopic equation is still a
 stochastic partial differential equation.
\end{remark}


\section{Long time effectivity  of the macroscopic model}\label{s5}

In this section we consider the long time effectivity  of the
homogenized system (\ref{limit}) in the autonomous case. It is
proved in section \ref{s4} that for fixed $T>0$ the macroscopic
behavior of the microscopic system (\ref{e1})-(\ref{e4}) can be
approximated by the macroscopic model (\ref{limit}) in the sense
of probability distribution. In fact we can show the long time
approximation. More specifically, we now prove that in the sense
of distribution, all  solutions of (\ref{e1})-(\ref{e4}) converge
to the unique stationary solution of (\ref{limit}) as
$T\rightarrow \infty$ and $\e\rightarrow 0$, under the assumption
that $f_\e\in H_\e$ and $g^i_\e\in V_\e$ are \emph{independent} of
time $t$ and
\begin{equation}\label{g11}
\sum_{i=1}^\infty |\nabla_{A_\e} g_\e^i(x)|^2_{\oplus_n H_\e}<C^*.
\end{equation}
Here $C^*$ is a positive constant independent of $\e$.

By the above assumptions, the property of  $a_{ij}$ and
$\beta_{ij}$, a standard argument (see \cite{PZ96}, Section 6)
yields that the system (\ref{abse}) and (\ref{limit}) have unique
stationary solutions $ u_\e^*(x, t)$ and
 $u^*(x,t)$ defined for $t>0$. We denote by  $\mu^*_\e$ and $\mu^*$ the
 distributions of $P_\e u^*_\e$ and $u^*$ in the space $H$, respectively. Then
  if $\mathbf{E}|u_\e^0|^2<\infty$ and $\mathbf{E}|u^0|^2<\infty$,
 \begin{equation}\label{conv0}
\Big|\int_H h d\mu_\e(t)-\int_H hd\mu^*_\e\Big|\leq
C(u_\e^0)e^{-\gamma t},\;\; t>0,
\end{equation}
\begin{equation}\label{conv}
\Big|\int_H h d\mu(t)-\int_H hd\mu^*\Big|\leq C(u^0)e^{-\gamma t},
\;\; t>0
\end{equation}
for some constant $\gamma>0$ and any $h: H\rightarrow \R^1$ with
 $\sup|h|\leq 1$ and $\textrm{Lip}(h)\leq 1$. Here $\mu_\e(t)=\mathcal{L}(P_\e u_\e(t, u_\e^0))$,
$\mu(t)=\mathcal{L}(u(t, \frac{u^0}{\vartheta}))$, and $C(u_\e^0)$
and $C(u^0)$ are positive constants   depending only on the initial
value $u_\e^0$ and $u^0$ respectively. The above convergence also
yields that $\mu_\e(t)$ and $\mu(t)$ weakly converges to $\mu^*_\e$
and $\mu^*$ respectively, as $t\rightarrow\infty $.

We will give some additional a \emph{priori} estimates which is
uniform with respect to $\e$ to ensure the tightness of the
stationary distributions.  For Banach space $U$ and $p>1$, we define
$W^{1,p}(0 ,T; U)$ as the space of functions $h\in L^p(0, T; U)$
such that
$$
 |h|^p_{W^{1, p}(0, T; U)}=|h|^p_{L^p(0, T; U)}
+\Big|\frac{dh}{dt}\Big|^p_{L^p(0, T; U)}<\infty.
$$
And for any $\alpha\in(0, 1)$, define $W^{\alpha, p}(0, T; U)$ as
the space of function $h\in L^p(0, T; U)$ such that
$$
 |h|^p_{W^{\alpha, p}(0, T; U)}=|h|^p_{L^p(0, T; U)}+
 \int_0^T\int_0^T\frac{|h(t)-h(s)|^p_U}{|t-s|^{1+\alpha
p}}dsdt<\infty.
$$
For $\rho\in(0,1)$, we denote by $C^\rho(0, T; U)$ the space of
functions  $h:[0,T]\rightarrow \mathcal{X}$ that are H\"{o}lder
continuous with exponent $\rho$.

\vskip 0.3cm

In the following part of this section we always assume that $f_\e$
and $g_\e^i$ are independent of time $t$ with (\ref{g11}) hold.
And for $T>0$  denote by $\mathfrak{u}^*_{\e, T}$ (respectively,
$\mathfrak{u}^*_T $) the distribution of stationary process $P_\e
u_\e^*(\cdot)$ (respectively, $u^*(\cdot)$) in the space $L^2(0,T;
H^1)$. Then we have the following result.

\begin{lemma} \label{longlimit}
For any $T>0$ the family $\mathfrak{u}^*_{\e, T}$ is tight in the
space $L^2(0, T; H^{2-\iota})$ with $\iota>0$.
\end{lemma}

\begin{proof}
 Since $u^*_\e$ is
stationary, by (\ref{est3}), we see that
\begin{equation}\label{1}
\mathbf{E}|u^*_\e|^2_{L^2(0, T; H_\e^2)}<C_T.
\end{equation}
Now represent $u^*_\e$ in the form
\begin{equation*}
u^*_\e(t)=u^*_\e(0)+\int_0^t\mathcal{A}_\e
u^*_\e(s)ds+\int_0^tf_\e(x)ds+\int_0^tg_\e(x)dW(s).
\end{equation*}
Also by the stationarity of  $u_\e^*$ and (\ref{est3})  we obtain
\begin{equation}\label{2}
\mathbf{E}\Big|\int_0^t\mathcal{A}_\e P_\e
u^*_\e(s)ds+\int_0^t\tilde{f}_\e(x)ds\Big|^2_{W^{1,2}(0, T;H)}\leq
C_T.
\end{equation}
Let $M_\e(s,t)=\int_s^t\tilde{g}_\e(x)dW(s)$. By Lemma 7.2 of
\cite{PZ92} and H$\ddot{o}$lder inequality, we derive that
\begin{eqnarray*}
\mathbf{E}|M_\e(s,t)|^4_{V_\e}\leq c
\Big(\int_s^t|\nabla_{A_\e}\tilde{g}_\e(x)|^2_{\mathcal{L}_2^Q}d\tau\Big)^2&\leq&
K(t-s)\int_s^t|\nabla_{A_\e}\tilde{g}_\e(x)|^4_{\mathcal{L}_2^Q}d\tau\\
&\leq& KC^{*2}|t-s|^2
\end{eqnarray*}
for $t\in [s,T]$, where $K$ is a positive constant independent of
$\e$, $s$ and $t$. Then
\begin{equation}\label{3}
\mathbf{E}\int_0^T|M_\e(0,t)|^4_{V_\e}dt\leq C_T
\end{equation}
and
\begin{eqnarray}\label{4}
\mathbf{E}\int_0^T\int_0^T\frac{|M_\e(0,t)-M_\e(0,s)|^4_{V_\e}}{|t-s|^{1+4\alpha}}dsdt\leq
C_T.
\end{eqnarray}

Combining (\ref{1})-(\ref{4}), and the compact embedding of
$$
L^2(0,T; H^2)\cap W^{1,2}(0,T; H)\subset L^2(0,T; H^{2-\iota})
$$
and
$$
L^2(0,T; H^2)\cap W^{\alpha, 4}(0,T; H^1)\subset L^2(0,T;
H^{2-\iota})
$$
we obtain the tightness of $\mathfrak{u}^*_{\e, T}$. This completes
the proof.
\end{proof}
The above lemma directly yields the following result
\begin{corollary}\label{longlimit1}
The family $\{\mu_\e^*\}$ is tight in the space $H^1$.
\end{corollary}

By Lemma \ref{longlimit}, for any fixed  $T>0$, the Skorohod
embedding theorem asserts that for any sequence $\{\e_n\}_n$ with
$\e_n\rightarrow 0$ as $n\rightarrow \infty$, there is subsequence
$\{\e_{n(k)}\}_k$, a new probability space $(\overline{\Omega},
\overline{\mathcal{F}}, \overline{\mathbb{P}})$ and random
variables $\overline{u}^*_{\e_{n(k)}}\in L^2(0, T; V_\e)$,
$\overline{u}^*\in L^2(0,T; H^1)$ such that
$$
 \mathcal{L}(P_\e \overline{u}^*_{\e_{n(k)}})=\mathfrak{u}^*_{\e_{n(k)}, T},\;\;
\mathcal{L}(\overline{u}^*)=\mathfrak{u}^*_T
$$
 and
 $$
\overline{u}^*_{\e_{n(k)}}\rightarrow \overline{u}^*,\;\;in\;\;
L^2(0,T; H^1)\;\; as \;\; k\rightarrow \infty.
 $$
Moreover $\overline{u}^*_{\e_{n(k)}}$ (respectively,
$\overline{u}^*$) is the unique stationary solution of  equation
(\ref{abse}) (respectively, (\ref{limit})) with $W$ replaced by
 $\overline{W}_k$ (respectively, $\overline{W}$). Here $\overline{W}_k$ and
 $\overline{W}$ are some Wiener processes defined on $(\overline{\Omega}, \overline{\mathcal{F}},
\overline{\mathbb{P}})$ with same distribution as $W$. Then by the
analysis of Section \ref{s4} and the uniqueness of the invariant
measure, we have
$$
\mathfrak{u}_{\e,T}^*\rightharpoonup \mathfrak{u}_T^*,\;\;
as\;\;\e\rightarrow 0
$$
for any $T>0$.

To show the long time effectivity,  let $u_\e(t)$, $t\geq 0$,  be a
weak  solution of system (\ref{e1})-(\ref{e4}) and define
$u_\e^t(\cdot)=u_\e(t+\cdot)$ which is in the space $L^2_{loc}(\R_+;
V_\e)$ by Theorem \ref{wellpose}. Then by (\ref{conv0})
$$
\mathcal{L}(P_\e u_\e^t(\cdot)) \rightharpoonup \mathcal{L}( P_\e
u_\e^*(\cdot)), \;\;t\rightarrow \infty
$$
in the space of probability measures on $L^2_{loc}(R_+; H^1)$.
Having the above analysis we draw the following result which implies
the long time effectivity of the homogenized effective equation
(\ref{limit}).

\begin{theorem}  (\textbf{Long time effectivity of macroscopic
model})\\
Assume that $f_\e\in H_\e$ and $g_\e^i\in V_\e$ are independent of
time $t$ with (\ref{g11}) hold, and further assume that (\ref{fe})
and (\ref{gee}) hold in $H$. Denote by $u_\e(t)$, $t\geq 0$,  the
solution of (\ref{e1})-(\ref{e4}) and $u^*$ the unique stationary
solution of (\ref{limit}). Then
\begin{equation}\label{long}
\lim_{\e\downarrow 0}\lim_{t\rightarrow \infty}\mathcal{L}(P_\e
u_\e^t(\cdot))=\mathcal{L}(u^*(\cdot)),
\end{equation}
where the limits are understood in the sense of weak convergence
of Borel probability measures in the space $L^2_{loc}(\R_+; H^1)$.
That is, the solution of (\ref{e1})-(\ref{e4}) converges to the
stationary solution of (\ref{limit}) in probability distribution
as $t \rightarrow \infty$ and
$\e\rightarrow 0$.\\
\end{theorem}

\begin{remark}
This theorem implies that the macroscopic  model (\ref{limit}) is
an effective approximation for the microscopic model (\ref{e1}),
on very long time scale. In other words, if we intend to
numerically simulate the long time behavior of the microscopic
model, we could just simulate the macroscopic  model as an
approximation when $\e$ is sufficiently small.

\end{remark}


\section{Effectivity in energy convergence   }\label{s6}

In the last two sections, we have considered finite time and long
time effectivity of the macroscopic model (\ref{limit}), in the
sense of convergence in probability distribution. In this section
we focus on the  finite time  effectivity of the macroscopic model
(\ref{limit}), but in the sense of   convergence in energy.
Namely, we show that the solution of   the microscopic model
(\ref{e1}) or (\ref{abse}),   converges   to the solution of the
macroscopic model (\ref{limit}), in an energy norm.

Let $u_\e$ be a weak solution of (\ref{abse}) and $u$ be a weak
solution of (\ref{limit}). We introduce the following energy
functionals:
\begin{equation}\label{Ee}
\mathcal{E}^\e(u_\e)(t)=\frac{1}{2}\mathbf{E}|\tilde{u}_\e|^2_H+\mathbf{E}\int_0^t\int_D\chi_{D_\e}A_\e\nabla
\big(P_\e u_\e(x,\tau)\big)\nabla \big(P_\e
u_\e(x,\tau)\big)dxd\tau
\end{equation}
and
\begin{equation}\label{E0}
\mathcal{E}^0(u)(t)=\frac{1}{2}\mathbf{E}|u|^2_H+\mathbf{E}\int_0^t\int_D
B\nabla u(x,\tau)\nabla u(x,\tau)dxd\tau.
\end{equation}
By the It$\acute{o}$ formula, it is clear that
$$
\mathcal{E}^\e(u_\e)(t)=\frac{1}{2}\mathbf{E}|\tilde{u}_\e^0|^2_H+\mathbf{E}\int_0^t\int_D\tilde{f}_\e(x,\tau)
\tilde{u}_\e(x,\tau)dxd\tau+\frac{1}{2}\mathbf{E}\int_0^t|\tilde{g}_\e(x,\tau)|^2_{\mathcal{L}_2^Q}d\tau
$$
and
$$
 \mathcal{E}^0(u)(t)=\frac{1}{2}\mathbf{E}|u^0|^2_H+\mathbf{E}\int_0^t\int_Df(x,\tau) u(x,\tau)dxd\tau
 +\frac{1}{2}\mathbf{E}\int_0^t|g(x,\tau)|^2_{\mathcal{L}_2^Q}d\tau.
$$

Then we have the following result on effectivity  of the
macroscopic model in the sense of  convergence in energy.

\begin{theorem}\label{energy}(\textbf{Effectivity in energy convergence})\\
Assume that (\ref{fe}) and (\ref{gee}) hold. If
$$
\tilde{u}^0_\e\rightarrow u^0,\;\;strongly \;in \; H, \text{as}\;
\e   \rightarrow 0,
$$
then
$$
\mathcal{E}^\e(u_\e)\rightarrow \mathcal{E}^0(u)\;\;\; in\;\;
C([0,T]), \text{as} \; \e \rightarrow 0.
$$
\end{theorem}



\begin{proof}
By the analysis of section \ref{s4}, for any $\delta>0$,
$u_\e\rightarrow u$ strongly in $L^2(0,T; H)$ on $\Omega_\delta$,
then by the arbitrariness of $\delta$, it is easy to see that
$$
\mathbf{E}\int_0^t\int_D\tilde{f}_\e(x,\tau)
\tilde{u}_\e(x,\tau)dxd\tau\rightarrow
\mathbf{E}\int_0^t\int_Df(x,\tau)
u(x,\tau)dxd\tau,\;\;for\;t\in[0, T].
$$

Then by $\tilde{g}_\e \rightharpoonup g$ weakly in
$L^2(0,t;\mathcal{L}_2^Q)$, we have
\begin{equation}\label{pointconv}
\mathcal{E}^\e(u_\e)(t)\rightarrow \mathcal{E}^0(u)(t)\;\; for\;\;
any\;\;t\in [0,T].
\end{equation}

We now only need to show that $\{\mathcal{E}^\e(u_\e)(t)\}_\e$ is
equicontinuous, as then the Ascoli-Arzela's theorem \cite{Dudley}
will imply the  result in the theorem.

In fact, given any $t\in[0,T]$, and $h>0$ small enough, we have

\begin{eqnarray*}
&& |\mathcal{E}^\e(u_\e)(t+h)-\mathcal{E}^\e(u_\e)(t)|\\
 &\leq& \Big|\mathbf{E}\int_t^{t+h} \int_D\tilde{f}_\e(x,\tau)\tilde{u}_\e(x,\tau)dxd\tau\Big|+
\mathbf{E}\int_t^{t+h}|\tilde{g}_\e(x,\tau)|^2_{\mathcal{L}_2^Q}d\tau\\
&\leq&\mathbf{E}\Big\{|\tilde{f}_\e|_{L^2(0,T;
H)}\int_t^{t+h}|\tilde{u}_\e(x,\tau)|^2_Hdxd\tau\Big\}+\mathbf{E}\int_t^{t+h}|\tilde{g}_\e(x,\tau)|^2_{\mathcal{L}_2^Q}d\tau.
\end{eqnarray*}

Noting that $\tilde{u}_\e\in L^2(0, T; H)$ a.s. and (\ref{g}), we
have
$$
|\mathcal{E}^\e(u_\e)(t+h)-\mathcal{E}^\e(u_\e)(t)|\rightarrow
0,\;\;as\; h\rightarrow 0,
$$
uniformly on $\e$, which means the equi-continuity  of the family
$\{\mathcal{E}^\e(u_\e)\}_\e$. This completes the proof.
\end{proof}


\section{Appendix: The homogenized matrix}
In this Appendix, we give the explicit expression of the
homogenized matrix $B$; for more details see \cite{CD99}. Let
$\chi^i$, $i=1,\cdots,n$ be the solutions of
\begin{eqnarray}
&&-\sum_{l,k=1}^n\frac{\partial}{\partial y_l}\Big(
a_{kl}\frac{\partial(\chi^i-y_i)}{\partial y_k}\Big)=0\;\;in\;\; Y^*
\label{asy1}
\end{eqnarray}
\begin{eqnarray}
&&\sum_{l,k=1}^n a_{kl}\frac{\partial(\chi^i-y_i)}{\partial y_k}
n_l=0\;\; on\;
\partial S\\
&&\chi^i \; \text{is}\;  Y-\text{periodic}. \label{asy3}
\end{eqnarray}
It is easy to calculate that $ \chi^i=-w_{e_i}+e_i$ with $\{e_i
\}_{i=1}^n$   the canonical basis of $\R^n$. Then
\begin{eqnarray}\label{beta}
\beta_{ij}=\frac{1}{|Y|}\int_Y\sum_{k=1}^na_{kj}\frac{\partial
w_{e_i}}{\partial y_k}dy=\frac{1}{|Y|}\int_Y
a_{ij}dy-\frac{1}{|Y|}\int_Y\sum_{k=1}^na_{kj}\frac{\partial
\chi_i}{\partial y_k}dy.
\end{eqnarray}
Moreover the operator $B=(\beta_{ij})$ satisfies the uniform
ellipticity condition: there is a constant $b>0$ such that
$$
\sum_{i,j=1}^n\beta_{ij}\xi_i\xi_j \geq b \sum_{i=1}^n\xi_i^2,
\;\; for\;\xi=(\xi_1,\cdots,\xi_n)\in\R^n.
$$


\vskip 0.51cm

 {\noindent{\bf{Acknowledgements:}} \vskip 0.25cm

The authors thank the referees for  very helpful suggestions and
comments.

\end{document}